# TIME SERIES ANALYSIS VIA MECHANISTIC MODELS

By Carles Bretó,[1] Daihai He, Edward L. Ionides[1,2,3]
and Aaron A. King[1,3]

*Universidad Carlos III de Madrid and University of Michigan*

The purpose of time series analysis via mechanistic models is to reconcile the known or hypothesized structure of a dynamical system with observations collected over time. We develop a framework for constructing nonlinear mechanistic models and carrying out inference. Our framework permits the consideration of *implicit* dynamic models, meaning statistical models for stochastic dynamical systems which are specified by a simulation algorithm to generate sample paths. Inference procedures that operate on implicit models are said to have the *plug-and-play* property. Our work builds on recently developed plug-and-play inference methodology for partially observed Markov models. We introduce a class of implicitly specified Markov chains with stochastic transition rates, and we demonstrate its applicability to open problems in statistical inference for biological systems. As one example, these models are shown to give a fresh perspective on measles transmission dynamics. As a second example, we present a mechanistic analysis of cholera incidence data, involving interaction between two competing strains of the pathogen *Vibrio cholerae*.

**1. Introduction.** A dynamical system is a process whose state varies with time. A mechanistic approach to understanding such a system is to write down equations, based on scientific understanding of the system, which describe how it evolves with time. Further equations describe the relationship of the state of the system to available observations on it. Mechanistic time series analysis concerns drawing inferences from the available data about

Received June 2008; revised August 2008.

[1]Supported in part by NSF Grant EF 0430120.

[2]Supported in part by NSF Grant DMS-08-05533.

[3]Supported in part by the RAPIDD program of the Science & Technology Directorate, Department of Homeland Security, and the Fogarty International Center, National Institutes of Health.

the hypothesized equations [Brillinger (2008)]. Questions of general interest include the following. Are the data consistent with a particular model? If so, for what range of values of model parameters? Does one mechanistic model describe the data better than another?

The defining principle of mechanistic modeling is that the model structure should be chosen based on scientific considerations, rather than statistical convenience. Although linear Gaussian models give an adequate representation of some processes [Durbin and Koopman (2001)], nonlinear behavior is an essential property of many systems. This leads to a need for statistical modeling and inference techniques applicable to rather general classes of processes. In the absence of alternative statistical methodology, a common approach to mechanistic investigations is to compare data, qualitatively or via some ad-hoc metric, with simulations from the model. It is a challenging problem, of broad scientific interest, to increase the range of mechanistic time series models for which formal statistical inferences, making efficient use of the data, can be made. Here, we develop a framework in which simulation of sample paths is employed as the basis for likelihood-based inference. Inferential techniques that require only simulation from the model (i.e., for which the model could be replaced by a black box which inputs parameters and outputs sample paths) have been called "equation free" [Kevrekidis, Gear and Hummer (2004), Xiu, Kevrekidis and Ghanem (2005)]. We will use the more descriptive expression "plug and play."

Plug-and-play inference techniques can be applied to any time series model for which a numerical procedure to generate sample paths is available. We call such models *implicit*, meaning that closed-form expressions for transition probabilities or sample paths are not required. The goal of this paper is to develop plug-and-play inference for a general class of implicitly specified stochastic dynamic models, and to show how this capability enables new and improved statistical analyses addressing current scientific debates. In other words, we introduce and demonstrate a framework for time series analysis via mechanistic models.

Here, we concern ourselves with partially observed, continuous-time, non-linear, Markovian stochastic dynamical systems. The particular combination of properties listed above is chosen because it arises naturally when constructing a mechanistic model. Although observations will typically be at discrete times, mechanistic equations describing underlying continuous time systems are most naturally described in continuous time. If all quantities important for the evolution of the system are explicitly modeled, then the future evolution of the system depends on the past only through the current state, that is, the system is Markovian. A stochastic model is prerequisite for mechanistic time series analysis, since chance variability is required to explain the difference between the data and the solution to noise-free deterministic equations. Statistical analysis is simpler if stochasticity can be confined to the observation process (the statistical problem becomes nonlinear regression) or if the stochastic dynamical system is perfectly observed [Basawa and Prakasa Rao (1980)]. Here we address the general case with both forms of stochasticity. Despite considerable work on such models [Anderson and Moore (1979), Liu (2001), Doucet, de Freitas and Gordon (2001), Cappé, Moulines and Rydén (2005)], statistical methodology which is readily applicable for a wide range of models has remained elusive. For example, Markov chain Monte Carlo and Monte Carlo Expectation-Maximization algorithms [Cappé, Moulines and Rydén (2005)] have technical difficulties handling continuous time dynamic models [Beskos et al. (2006)]; these two approaches also lack the plug-and-play property.

Several inference techniques have previously been proposed which are compatible with plug-and-play inference from partially observed Markov processes. Nonlinear forecasting [Kendall et al. (1999)] is a method of simulated moments which approximates the likelihood. Iterated filtering is a recently developed method [Ionides, Bretó and King (2006)] which provides a way to calculate a maximum likelihood estimate via sequential Monte Carlo, a plug and play filtering technique. Approximate Bayesian sequential Monte Carlo plug-and-play methodologies [Liu and West (2001), Toni et al. (2008)] have also been proposed.

In Section 2 we introduce a new and general class of implicitly specified models. Section 3 is concerned with inference methodology and includes a review of the iterated filtering approach of Ionides, Bretó and King (2006). Section 4 discusses the role of our modeling and inference framework for the analysis of biological systems. Two concrete examples are developed, investigating measles (Section 4.1) and cholera (Section 4.2). Section 5 is a concluding discussion. The motivating examples in this paper have led to an emphasis on modeling infectious diseases. However, the issue of mechanistic modeling of time series data occurs in many other contexts. Indeed, it is too widespread to give a comprehensive review and we instead list some examples: molecular biochemistry [Kou, Xie and Liu (2005)]; wildlife ecology [Newman and Lindley (2006)]; cell biology [Ionides et al. (2004)]; economics [Fernández-Villaverde and Rubio-Ramírez (2005)]; signal processing [Arulampalam et al. (2002)]; data assimilation for numerical models [Houtekamer and Mitchell (2001)]. The study of infectious disease, however, has a long history of motivating new modeling and data analysis methodology [Kermack and McKendrick (1927), Bartlett (1960), Anderson and May (1991), Finkenstädt and Grenfell (2000), Ionides, Bretó and King (2006), King et al. (2008b)]. The freedom to carry out formal statistical analysis based on mechanistically motivated, nonlinear, nonstationary, continuous time stochastic models is a new development which promises to be a useful tool for a variety of applications.

**2. Compartment models with stochastic rates.** Many mechanistic models can be viewed in terms of flows between compartments [Jacquez (1996), Matis and Kiffe (2000)]. Here, we introduce a class of implicitly specified stochastic compartment models; widespread biological applications of these models will be discussed in Section 4, with broader relevance and further generalizations discussed in Section 5. The reader may choose initially to pass superficially through the technical details of this section. We present a general model framework which is, at once, an example of an implicitly specified mechanistic model, a necessary prelude to our following data analyses, and a novel class of Markov processes requiring some formal mathematical treatment.

A general compartment model is a vector-valued process $X(t) = (X_1(t), \ldots, X_c(t))$ denoting the (integer or real-valued) counts in each of $c$ compartments. The basic characteristic of a compartment model is that $X(t)$ can be written in terms of the flows $N_{ij}(t)$ from $i$ to $j$, via a "conservation of mass" identity:

$$\text{(1)} \qquad X_i(t) = X_i(0) + \sum_{j \neq i} N_{ji}(t) - \sum_{j \neq i} N_{ij}(t).$$

Each *flow* $N_{ij}$ is associated with a *rate* function $\mu_{ij} = \mu_{ij}(t, X(t))$. There are many ways to develop concrete interpretations of such a compartment model. For the remainder of this section, we take $X_i(t)$ to be nonnegative integer-valued, so $X(t)$ models a population divided into $c$ disjoint categories and $\mu_{ij}$ is the rate at which each individual in compartment $i$ moves to $j$. In this context, it is natural to require that $\{N_{ij}(t), 1 \leq i \leq c, 1 \leq j \leq c\}$ is a collection of nondecreasing integer-valued stochastic processes satisfying the constraint $X_i(t) \geq 0$ for all $i$ and $t$. The conservation equation (1) makes the compartment model closed in the sense that individuals cannot enter or leave the population. However, processes such as immigration, birth or death can be modeled via the introduction of additional source and sink compartments.

We wish to introduce white noise to model stochastic variation in the rates (discussion of this decision is postponed to Sections 4 and 5). We refer to *white noise* as the derivative of an *integrated noise* process with stationary independent increments [Karlin and Taylor (1981)]. The integral of a white noise process over an interval is thus well defined, even when the sample paths of the integrated noise process are not formally differentiable. Specifically, we introduce a collection of integrated noise processes $\{\Gamma_{ij}(t), 1 \leq i \leq c, 1 \leq j \leq c\}$ with the properties:

(P1) Independent increments: The collection of increments $\{\Gamma_{ij}(t_2) - \Gamma_{ij}(t_1), 1 \leq i \leq c, 1 \leq j \leq c\}$ is presumed to be independent of $\{\Gamma_{ij}(t_4) - \Gamma_{ij}(t_3), 1 \leq i \leq c, 1 \leq j \leq c\}$ for all $t_1 < t_2 < t_3 < t_4$.

1. Divide the interval $[0, T]$ into $N$ intervals of width $\delta = T/N$
2. Set initial value $X(0)$
3. FOR $n = 0$ to $N - 1$
4.     Generate noise increments $\{\Delta\Gamma_{ij} = \Gamma_{ij}(n\delta + \delta) - \Gamma_{ij}(n\delta)\}$
5.     Generate process increments
   $$(\Delta N_{i1}, \ldots, \Delta N_{i,i-1}, \Delta N_{i,i+1}, \Delta N_{ic}, R_i)$$
   $$\sim \text{Multinomial}(X_i(n\delta), p_{i1}, \ldots, p_{i,i-1},$$
   $$p_{i,i+1}, \ldots, p_{ic}, 1 - \textstyle\sum_{k \neq i} p_{ik})$$
   where $p_{ij} = p_{ij}(\{\mu_{ij}(n\delta, X(n\delta))\}, \{\Delta\Gamma_{ij}\})$ is given in (3)
6.     Set $X_i(n\delta + \delta) = R_i + \sum_{j \neq i} \Delta N_{ji}$
7. END FOR

FIG. 1. *Euler scheme for a numerical solution of the Markov chain specified by (2). In steps 5 and 6, $R_i$ counts the individuals who remain in compartment $i$ during the current Euler increment.*

(P2) Stationary increments: The collection of increments $\{\Gamma_{ij}(t_2) - \Gamma_{ij}(t_1), 1 \leq i \leq c, 1 \leq j \leq c\}$ has a joint distribution depending only on $t_2 - t_1$.

(P3) Nonnegative increments: $\Gamma_{ij}(t_2) - \Gamma_{ij}(t_1) \geq 0$ for $t_2 > t_1$.

We have not assumed that different integrated noise processes $\Gamma_{ij}$ and $\Gamma_{kl}$ are independent; their increments could be correlated, or even equal. These integrated noise processes define a collection of noise processes given by $\xi_{ij}(t) = \frac{d}{dt}\Gamma_{ij}(t)$. Since $\Gamma_{ij}(t)$ is increasing, $\xi_{ij}(t)$ is nonnegative and $\mu_{ij}\xi_{ij}(t)$ can be interpreted as a rate with multiplicative white noise. In this context, it is natural to assume the following:

(P4) Unbiased multiplicative noise: $E[\Gamma_{ij}(t)] = t$.

At times, we may further assume one or more of the following properties:

(P5) Partially independent noises: For each $i$, $\{\Gamma_{ij}(t)\}$ is independent of $\{\Gamma_{ik}(t)\}$ for all $j \neq k$.

(P6) Independent noises: $\{\Gamma_{ij}(t)\}$ is independent of $\{\Gamma_{kl}(t)\}$ for all pairs $(i,j) \neq (k,l)$.

(P7) Gamma noises: Marginally $\Gamma_{ij}(t + \delta) - \Gamma_{ij}(t) \sim \text{Gamma}(\delta/\sigma_{ij}^2, \sigma_{ij}^2)$, the gamma distribution whose shape parameter is $\delta/\sigma_{ij}^2$ and scale parameter is $\sigma_{ij}^2$, with corresponding mean $\delta$ and variance $\delta\sigma_{ij}^2$. We call $\sigma_{ij}^2$ an *infinitesimal variance* parameter [Karlin and Taylor (1981)].

The choice of gamma noise in (P7) gives a convenient concrete example. A wide range of Lévy processes [Sato (1999)] could be alternatively employed.

We proceed to construct a compartment model as a continuous time Markov chain via the limit of coupled discrete-time multinomial processes with random rates. Similar Euler multinomial schemes (without noise in the rate function) are a standard numerical approach for studying population dynamics [Cai and Xu (2007)]. The representation of our model given in (2) is *implicit* since numerical solution is available to arbitrary precision via evaluating the coupled multinomial processes in a discrete time-step Euler scheme (described in Figure 1). Let $\Delta N_{ij} = N_{ij}(t + \delta) - N_{ij}(t)$ and $\Delta\Gamma_{ij} = \Gamma_{ij}(t + \delta) - \Gamma_{ij}(t)$. We suppose that

$$P[\Delta N_{ij} = n_{ij}, \text{for all } 1 \leq i \leq c, 1 \leq j \leq c, i \neq j \mid X(t) = (x_1, \ldots, x_c)]$$

$$(2) \qquad = E\left[\prod_{i=1}^{c}\left\{\binom{x_i}{n_{i1}\cdots n_{ii-1}n_{ii+1}\cdots n_{ic}r_i}\left(1 - \sum_{k\neq i}p_{ik}\right)^{r_i}\prod_{j\neq i}p_{ij}^{n_{ij}}\right\}\right]$$

$$\qquad + o(\delta),$$

where $r_i = x_i - \sum_{k\neq i} n_{ik}$, $\binom{n}{n_1\cdots n_c}$ is a multinomial coefficient and

$$(3) \qquad p_{ij} = p_{ij}(\{\mu_{ij}(t,x)\}, \{\Delta\Gamma_{ij}(t)\})$$

$$\qquad = \left(1 - \exp\left\{-\sum_k \mu_{ik}\Delta\Gamma_{ik}\right\}\right)\mu_{ij}\Delta\Gamma_{ij}\Big/\sum_k \mu_{ik}\Delta\Gamma_{ik},$$

with $\mu_{ij} = \mu_{ij}(t, x)$. Theorem A.1, which is stated in Appendix A and proved in a supplement to this article [Bretó et al. (2009)], shows that (2) defines a continuous time Markov chain when the conditions (P1)–(P5) hold. A finite-state continuous time Markov chain is specified by its infinitesimal transition probabilities [Brémaud (1999)], which are in turn specified by (2). Theorem A.2, also stated in Appendix A and proved in the supplement [Bretó et al. (2009)], determines the infinitesimal transition probabilities resulting from (2) supposing the conditions (P1)–(P7). When the infinitesimal transition probabilities can be calculated exactly, exact simulation methods are available [Gillespie (1977)]. In practice, numerical schemes based on Euler approximations may be preferable—Euler schemes for Markov chain compartment models have been proposed based on Poisson [Gillespie (2001)], binomial [Tian and Burrage (2004)] and multinomial [Cai and Xu (2007)] approximations. Our choice of a model for which convenient numerical solutions are available (e.g., via the procedure in Figure 1) comes at the expense of difficulty in computing analytic properties of the implicitly-defined continuous-time process. However, since the properties of the model will be investigated by simulation, via a plug-and-play methodology, the analytic properties of the continuous-time process are of relatively little interest.

For the gamma noise in (P7), the special case where $\sigma_{ij} = 0$ is taken to correspond to $\xi_{ij}(t) = 1$. If $\sigma_{ij} = 0$ for all $i$ and $j$, then (2) becomes the Poisson system widely used to model demographic stochasticity in population models [Brémaud (1999), Bartlett (1960)]. We therefore call a process defined by (2) a Poisson system with stochastic rates. Constructions similar to Theorem A.1 are standard for Poisson systems [Brémaud (1999)], but here care is required to deal with the novel inclusion of white noise in the rate process. Our formulation for adding noise to Poisson systems can be seen as a generalization of subordinated Lévy processes [Sato (1999)], though we are not aware of previous work on the more general Markov processes constructed here. It is only the recent development of plug-and-play inference methodology that has led to the need for flexible Markov chain models with random rates.

2.1. *Comments on the role of numerical solutions based on discretizations.* In this section we have proposed employing a discrete-time approximation to a continuous-time stochastic process. Numerical solutions based on discretizations of space and time are ubiquitous in the applied mathematical sciences and engineering. A standard technique is to investigate whether further reduction in the size of the discretization substantially affects the conclusions of the analysis. When sufficiently fine discretization is not computationally feasible, the numerical solutions may still have some value. Climate modeling and numerical weather prediction are examples of this: such systems have important dynamic behavior at scales finer than any feasible discretization, but numerical models nevertheless have a scientific role to play [Solomon et al. (2007)].

When numerical modeling is used as a scientific tool, conclusions about the limiting continuous time model will be claimed based on properties of the model that are determined by simulation of realizations from the discretized model. Such conclusions depend on the assumption that properties of the numerical solution which are stable as the numerical approximation timestep, $\delta$, approaches 0 should indeed be properties of the limiting continuous time process. This need not always be true, which is one reason why analytic properties, such as Theorems A.1 and A.2, are valuable.

From another point of view, an argument for being content with a numerical approximation to (2) for sufficiently small $\delta$ is that there may be no scientific reason to prefer a true continuous time model over a fine discretization. For example, when modeling year-to-year population dynamics, continuous time models of adequate simplicity for data analysis typically will not include diurnal effects. Thus, there is no particular reason to think the continuous time model more credible than a discrete time model with a step of one day. One can think of a set of equations defining a continuous time process, combined with a specified discretization, as a way of writing down a discrete time model, rather than treating the continuous time model as a gold standard against which all discretizations must be judged.

**3. Plug-and-play inference methodology.** We suppose now that the dynamical system depends on some unknown parameter vector $\theta \in \mathbb{R}^{d_\theta}$, so that $\mu_{ij} = \mu_{ij}(t, x, \theta)$ and $\sigma_{ij} = \sigma_{ij}(t, x, \theta)$. Inference on $\theta$ is to be made based on observations $y_{1:N} = (y_1, \ldots, y_N)$ made at times $t_{1:N} = (t_1, \ldots, t_N)$, with $y_n \in \mathbb{R}^{d_y}$. Conditionally on $X(t_1), \ldots, X(t_N)$, we suppose that the observations are drawn independently from a density $g(y_n | X(t_n), \theta)$. Likelihood-based inference can be carried out for the framework of Section 2 using the iterated filtering methodology proposed by Ionides, Bretó and King (2006), implemented as described in Figure 2. Iterated filtering is a technique to maximize the likelihood for a partially observed Markov model, permitting calculation of maximum likelihood point estimates, confidence intervals (via profile likelihood, bootstrap or Fisher information) and likelihood ratio hypothesis tests. Iterated filtering has been developed in response to challenges arising in ecological and epidemiological data analysis [Ionides, Bretó and King (2006, 2008), King et al. (2008b)], and appears here for the first time in the statistical literature. We refer to Ionides, Bretó and King (2006) for the mathematical results concerning the iterated filtering algorithm in Figure 2. We proceed to review the methodology and its heuristic motivation, to discuss implementation issues, and to place iterated filtering in the context of alternative statistical methodologies.

For nonlinear non-Gaussian partially observed Markov models, the likelihood function can typically be evaluated only inexactly and at considerable computational expense. The iterated filtering procedure takes advantage of the partially observed Markov structure to enable computationally efficient maximization. A useful property of partially observed Markov models is that, if the parameter $\theta$ is replaced by a random walk $\theta_{1:N}$, with $E[\theta_0] = \theta$ and $E[\theta_n | \theta_{n-1}] = \theta_{n-1}$ for $n > 1$, the calculation of $\hat{\theta}_n = E[\theta_n | y_{1:n}]$ and $V_n = \mathrm{Var}(\theta_n | y_{1:n-1})$ is a well-studied and computationally convenient filtering problem [Kitagawa (1998), Liu and West (2001)]. Additional stochasticity of this kind is introduced in steps 4 and 12 of Figure 2. This leads to time-varying parameter estimates, so $\bar{\theta}_i(t_n)$ in Figure 2 is an estimate of $\theta_i$ depending primarily on the data at and shortly before time $t_n$. The updating rule in step 16, giving an appropriate way to combine these temporally local estimates, is the main innovative component of the procedure. Ionides, Bretó and King (2006) showed that this algorithm converges to the maximum of the likelihood function, under sufficient regularity conditions to justify a Taylor series expansion argument. Only the mean and variance of the stochasticity added in steps 4 and 12 play a role in the limit as $n$ increases. The specific choice of the normal distribution for steps 4 and 12 of Figure 2 is therefore unimportant, but does require that the parameter space is unbounded. This is achieved by reparameterizing where necessary; we use a log transform for positive parameters and a logit transform for parameters lying in the interval $[0, 1]$.

Steps 2, 11 and 17 of Figure 2 concern the initial values of the state variables. For stationary processes, one can think of these as unobserved random

MODEL INPUT: $f(\cdot)$, $g(\cdot|\cdot)$, $y_1, \ldots, y_N$, $t_0, \ldots, t_N$

ALGORITHMIC PARAMETERS: integers $J$, $L$, $M$;
scalars $0 < a < 1$, $b > 0$; vectors $X_I^{(1)}$, $\theta^{(1)}$;
positive definite symmetric matrices $\Sigma_I$, $\Sigma_\theta$.

1. FOR $m = 1$ to $M$
2.      $X_I(t_0, j) \sim N[X_I^{(m)}, a^{m-1}\Sigma_I]$,    $j = 1, \ldots, J$
3.      $X_F(t_0, j) = X_I(t_0, j)$
4.      $\theta(t_0, j) \sim N[\theta^{(m)}, ba^{m-1}\Sigma_\theta]$
5.      $\bar{\theta}(t_0) = \theta^{(m)}$
6.      FOR $n = 1$ to $N$
7.          $X_P(t_n, j) = f(X_F(t_{n-1}, j), t_{n-1}, t_n, \theta(t_{n-1}, j), W)$
8.          $w(n, j) = g(y_n | X_P(t_n, j), t_n, \theta(t_{n-1}, j))$
9.          draw $k_1, \ldots, k_J$ such that
             $\mathrm{Prob}(k_j = i) = w(n, i) / \sum_\ell w(n, \ell)$
10.     $X_F(t_n, j) = X_P(t_n, k_j)$
11.     $X_I(t_n, j) = X_I(t_{n-1}, k_j)$
12.     $\theta(t_n, j) \sim N[\theta(t_{n-1}, k_j), a^{m-1}(t_n - t_{n-1})\Sigma_\theta]$
13.     Set $\bar{\theta}_i(t_n)$ to be the sample mean of $\{\theta_i(t_{n-1}, k_j), j = 1, \ldots, J\}$
14.     Set $V_i(t_n)$ to be the sample variance of $\{\theta_i(t_n, j), j = 1, \ldots, J\}$
15.    END FOR
16.    $\theta_i^{(m+1)} = \theta_i^{(m)} + V_i(t_1) \sum_{n=1}^N V_i^{-1}(t_n)(\bar{\theta}_i(t_n) - \bar{\theta}_i(t_{n-1}))$
17.    Set $X_I^{(m+1)}$ to be the sample mean of $\{X_I(t_L, j), j = 1, \ldots, J\}$
18. END FOR

RETURN
maximum likelihood estimate for parameters, $\hat{\theta} = \theta^{(M+1)}$
maximum likelihood estimate for initial values, $\hat{X}(t_0) = X_I^{(M+1)}$
maximized conditional log likelihood estimates, $\ell_n(\hat{\theta}) = \log(\sum_j w(n, j)/J)$
maximized log likelihood estimate, $\ell(\hat{\theta}) = \sum_n \ell_n(\hat{\theta})$

Fig. 2. *Implementation of likelihood maximization by iterated filtering. $N[\mu, \Sigma]$ corresponds to a normal random variable with mean vector $\mu$ and covariance matrix $\Sigma$; $X(t_n)$ takes values in $\mathbb{R}^{d_x}$; $y_n$ takes values in $\mathbb{R}^{d_y}$; $\theta$ takes values in $\mathbb{R}^{d_\theta}$ and has components $\{\theta_i, i = 1, \ldots, d_\theta\}$; $f(\cdot)$ is the transition rule described in (4); $g(\cdot|\cdot)$ is the measurement density for the observations $y_{1:N}$.*

variables drawn from the stationary distribution. However, for nonstationary processes (such as those considered in Sections 4.1 and 4.2, and any process modeled conditional on measured covariates) these initial values are treated as unknown parameters. These parameters require special attention, despite not usually being quantities of primary scientific interest, since the information about them is concentrated at the beginning of the time series, whereas the computational benefit of iterated filtering arises from combining information accrued through time. Steps 2, 11 and 17 implement a fixed lag smoother [Anderson and Moore (1979)] to iteratively update the initial value estimates. The value of the fixed lag (denoted by $L$ in Figure 2) should be chosen so that there is negligible additional information about the initial values after time $t_L$. Choosing $L$ too large results in slower convergence, choosing $L$ too small results in bias.

Iterated filtering, characterized by the updating rule in step 16 of Figure 2, can be implemented via any filtering method. The procedure in Figure 2 employs a basic sequential Monte Carlo filter which we found to be adequate for the examples in Section 4 and also for previous data analyses [Ionides, Bretó and King (2006), King et al. (2008b)]. Many extensions and generalizations of sequential Monte Carlo have been proposed [Arulampalam et al. (2002), Doucet, de Freitas and Gordon (2001), Del Moral, Doucet and Jasra (2006)] and could be employed in an iterated filtering algorithm. If the filtering technique is plug-and-play, then likelihood maximization by iterated filtering also has this property. Basic sequential Monte Carlo filtering techniques do have the plug-and-play property, since only simulations from the transition density of the dynamical system are required and not evaluation of the density itself. Although sequential Monte Carlo algorithms are usually written in terms of transition densities [Arulampalam et al. (2002), Doucet, de Freitas and Gordon (2001)], we emphasize the plug-and-play property of the procedure in Figure 2 by specifying a Markov process at a sequence of times $t_0 < t_1 < \cdots < t_N$ via a recursive transition rule,

$$(4) \qquad X(t_n) = f(X(t_{n-1}), t_{n-1}, t_n, \theta, W).$$

Here, it is understood that $W$ is some random variable which is drawn independently each time $f(\cdot)$ is evaluated. In the context of the plug-and-play philosophy, $f(\cdot)$ is the algorithm to generate a simulated sample path of $X(t)$ at the discrete times $t_1, \ldots, t_N$ given an initial value $X(t_0)$.

To check whether global maximization has been achieved, one can and should consider various different starting values [i.e., $\theta^{(1)}$ and $X_I^{(1)}$ in Figure 2]. Attainment of a local maximum can be checked by investigation of the likelihood surface local to an estimate $\hat{\theta}$. Such an investigation can also give rise to standard errors, and we describe here how this was carried out for the results in Table 2. We write $\ell(\theta)$ for the log likelihood function, and we call a graph of $\ell(\hat{\theta} + z\delta_i)$ against $\hat{\theta}_i + z$ a sliced likelihood plot; here, $\delta_i$ is a vector of zeros with a one in the $i$th position, and $\theta$ has components $\{\theta_1, \ldots, \theta_{d_\theta}\}$. If $\hat{\theta}$ is located at a local maximum of each sliced likelihood, then $\hat{\theta}$ is a local maximum of $\ell(\theta)$, supposing $\ell(\theta)$ is continuously differentiable. We check this by evaluating $\ell(\hat{\theta} + z_{ij}\delta_i)$ for a collection $\{z_{ij}\}$ defining a neighborhood of $\hat{\theta}$. The likelihood is evaluated with Monte Carlo error, as described in Figure 2, with $\Sigma_I = 0$, $\Sigma_\theta = 0$ and $M = 1$. Therefore, it is necessary to make a smooth approximation to the sliced likelihood [Ionides (2005)] based on the available evaluations. The size of the neighborhood (specified by $\{z_{ij}\}$) and the size of the Monte Carlo sample (specified by $J$ in Figure 2) should be large enough that the local maximum for each slice is clearly identified. Computing sliced likelihoods requires moderate computational effort, linear in the dimension of $\theta$. As a by-product of the sliced likelihood computation, one has access to the conditional log likelihood values, defined in Figure 2 and written here as $\ell_{nij} = \ell_n(\hat{\theta} + z_{ij}\delta_i)$. Regressing $\ell_{nij}$ on $z_{ij}$ for each fixed $i$ and $n$ gives rise to estimates $\hat{\ell}_{ni}$ for the partial derivatives of the conditional log likelihoods. Standard errors of parameters are found from the estimated observed Fisher information matrix [Barndorff-Nielsen and Cox (1994)], with entries given by $\hat{\mathcal{I}}_{ik} = \sum_n \dot{\ell}_{ni}\dot{\ell}_{nk}$. We prefer profile likelihood calculations, such as Figure 7, to derive confidence intervals for quantities of particular interest. However, standard errors derived from estimating the observed Fisher information involve substantially less computation.

Parameter estimation for partially observed nonlinear Markov processes has long been a challenging problem, and it is premature to expect a fully automated statistical procedure. The implementation of iterated filtering in Figure 2 employs algorithmic parameters which require some trial and error to select. However, once the likelihood has been demonstrated to be successfully maximized, the algorithmic parameters play no role in the scientific interpretation of the results.

Other plug-and-play inference methodologies applicable to the models of Section 2 have been developed. Nonlinear forecasting [Kendall et al. (1999)] has neither the statistical efficiency of a likelihood-based method nor the computational efficiency of a filtering-based method. The Bayesian sequential Monte Carlo approximation of Liu and West (2001) combines likelihood-based inference with a filtering algorithm, but is not supported by theoretical guarantees comparable to those presented by Ionides, Bretó and King (2006) for iterated filtering. A recently developed plug-and-play approach to approximate Bayesian inference [Sisson, Fan and Tanaka (2007)] has been applied to partially observed Markov processes [Toni et al. (2008)]. Other recent developments in Bayesian methodology for partially observed Markov processes include Newman et al. (2008), Cauchemez and Ferguson (2008), Cauchemez et al. (2008), Beskos et al. (2006), Polson, Stroud and Muller (2008), Boys, Wilkinson and Kirkwood (2008). This research has been motivated by the inapplicability of general Bayesian software, such as WinBUGS [Lunn et al. (2000)], for many practical inference situations [Newman et al. (2008)].

A numerical implementation of iterated filtering is available via the software package `pomp` [King, Ionides and Bretó (2008a)] which operates in the free, open-source, R computing environment [R Development Core Team (2006)]. This package contains a tutorial vignette as well as further examples of mechanistic time series models. The data analyses of Section 4 were carried out using `pomp`, in which the algorithms in Figures 1 and 2 are implemented via the functions `reulermultinom` and `mif` respectively.

**4. Time series analysis for biological systems.** Mathematical models for the temporal dynamics of biological populations have long played a role in understanding fluctuations in population abundance and interactions between species [Bjornstad and Grenfell (2001), May (2004)]. When using models to examine the strength of evidence concerning rival hypotheses about a system, a model is typically required to capture not just the qualitative features of the dynamics but also to explain quantitatively all the available observations on the system. A critical aspect of capturing the statistical behavior of data is an adequate representation of stochastic variation, which is a ubiquitous component of biological systems. Stochasticity can also play an important role in the qualitative dynamic behavior of biological systems [Coulson, Rohani and Pascual (2004), Alonso, McKane and Pascual (2007)]. Unpredictable event times of births, deaths and interactions between individuals result in random variability known as *demographic stochasticity* (from a microbiological perspective, the individuals in question might be cells or large organic molecules). The environmental conditions in which the system operates will fluctuate considerably in all but the best experimentally controlled situations, resulting in *environmental stochasticity*. The framework of Section 2 provides a general way to build the phenomenon of environmental stochasticity into continuous-time population models, via the inclusion of variability in the rates at which population processes occur. To our knowledge, this is the first general framework for continuous time, discrete population dynamics which allows for both demographic stochasticity [infinitesimal variance equal to the infinitesimal mean; see Karlin and Taylor (1981) and Appendix B] and environmental stochasticity [infinitesimal variance greater than the infinitesimal mean; see supporting online material Bretó et al. (2009)].

From the point of view of statistical analysis, environmental stochasticity plays a comparable role for dynamic population models to the role played by over-dispersion in generalized linear models. Models which do not permit consideration of environmental stochasticity lead to strong assumptions about the levels of stochasticity in the system. This relationship is discussed further in Section 5. For generalized linear models, over-dispersion is commonplace, and failure to account properly for it can give rise to misleading conclusions [McCullagh and Nelder (1989)]. Phrased another way, including sufficient stochasticity in a model to match the unpredictability of the data is essential if the model is to be used for forecasting, or predicting a quantitative range of likely effects of an intervention, or estimating unobserved components of the system.

We present two examples. First, Section 4.1 demonstrates the role of environmental stochasticity in measles transmission dynamics, an extensively studied and relatively simple biological system. Second, Section 4.2 analyzes data on competing strains of cholera to demonstrate the modeling framework of Section 2 and the inference methodology of Section 3 on a more complex system.

4.1. *Environmental stochasticity in measles epidemics.* The challenges of moving from mathematical models, which provide some insight into the system dynamics, to statistical models, which both capture the mechanistic basis of the system and statistically describe the data, are well documented by a sequence of work on the dynamics of measles epidemics [Bartlett (1960), Anderson and May (1991), Finkenstädt and Grenfell (2000), Bjornstad, Finkenstadt and Grenfell (2002), Morton and Finkenstadt (2005), Cauchemez and Ferguson (2008)]. Measles is no longer a major developed world health issue but still causes substantial morbidity and mortality, particularly in sub-Saharan Africa [Grais et al. (2006), Conlan and Grenfell (2007)]. The availability of excellent data before the introduction of widespread vaccination has made measles a model epidemic system. Recent attempts to analyze population-level time series data on measles epidemics via mechanistic dynamic models have, through statistical expediency, been compelled to use a discrete-time dynamic model using timesteps synchronous with the reporting intervals [Finkenstädt and Grenfell (2000), Bjornstad, Finkenstadt and Grenfell (2002), Morton and Finkenstadt (2005)]. Such discrete time models risk incorporating undesired artifacts [Glass, Xia and Grenfell (2003)]. The first likelihood-based analysis via continuous time mechanistic models, incorporating only demographic stochasticity, was published while this paper was under review [Cauchemez and Ferguson (2008)]. From another perspective, the properties of stochastic dynamic epidemic models have been studied extensively in the context of continuous time models with only demographic stochasticity [Bauch and Earn (2003), Dushoff et al. (2004), Wearing, Rohani and Keeling (2005)]. We go beyond previous approaches, by demonstrating the possibility of carrying out modeling and data analysis via continuous time mechanistic models with both demographic and environmental stochasticity. For comparison with the work of Cauchemez and Ferguson (2008), we analyze measles epidemics occurring in London, England during the pre-vaccination era. The data, reported cases from 1948 to 1964, are shown in Figure 3.

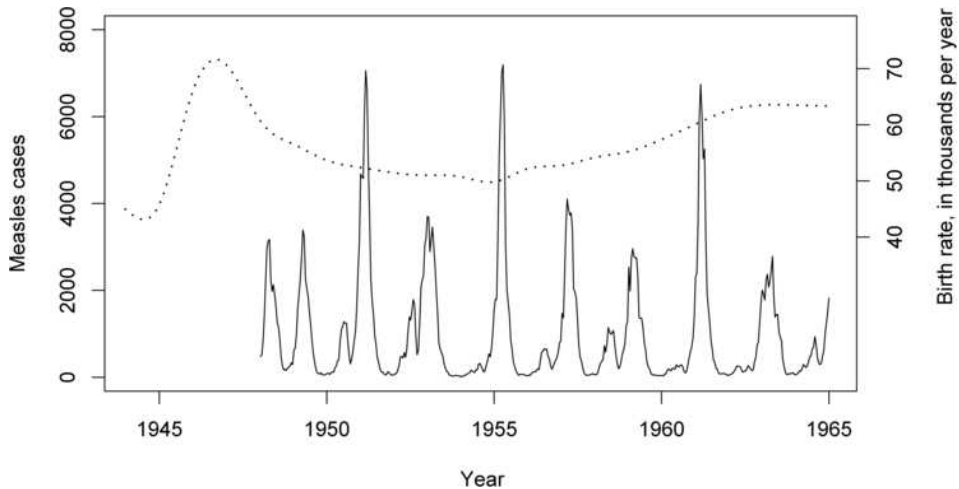

FIG. 3.  *Biweekly recorded measles cases (solid line) and birth rate (dotted line, estimated by smooth interpolation of annual birth statistics) for London, England.*

Measles has a relatively simple natural history, being a highly infectious disease with lifelong immunity following infection. As such, it was historically a childhood disease, with transmission occurring primarily in school environments. A basic model for measles has children becoming susceptible to infection upon reaching school age [here taken to be $\tau = 4$ year, following Cauchemez and Ferguson (2008)]. This susceptible group is described by a compartment $S$ containing, at time $t$, a number of individuals $S(t)$. Upon exposure to infection, a transition occurs to compartment $E$ in which individuals are infected but not yet infectious. The disease then progresses to an infectious state $I$. Individuals are finally removed to a state $R$, due to bed-rest and subsequent recovery. This sequence of compartments is displayed graphically in Figure 4. We proceed to represent this system as a Markov chain with stochastic rates, via the notation of Section 2, with $X(t) = (S(t), E(t), I(t), R(t), B(t), D(t))$. Compartments $B$ and $D$ are introduced for demographic considerations: births are represented by transitions from $B$ to $S$, with an appropriate delay $\tau$, and deaths by transitions into $D$.

The seasonality of the transmissibility has been found to be well described by whether or not children are in school [Fine and Clarkson (1982),

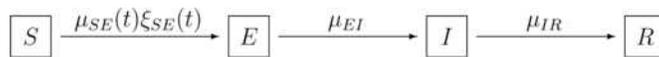

FIG. 4.  *Flow diagram for measles. Each individual host falls in one compartment: $S$, susceptible; $E$, exposed and infected but not yet infectious; $I$, infectious; $R$, removed and subsequently recovered and immune. Births enter $S$ after a delay $\tau$, and all individuals have a mortality rate $m$.*

Bauch and Earn (2003)]. Thus, we define a transmissibility function $\beta(t)$ by

$$(5) \qquad \beta(t) = \begin{cases} \beta_H : (t) = 7\text{--}99,\ 116\text{--}199,\ 252\text{--}299,\ 308\text{--}355, \\ \beta_L : d(t) = 356\text{--}6,\ 100\text{--}115,\ 200\text{--}251,\ 300\text{--}307, \end{cases}$$

where $d(t)$ is the integer-valued day (1–365) corresponding to the real-valued time $t$ measured in years. This functional form allows reduced transmission during the Christmas vacation (days 356–365 and 1–6), Easter vacation (100–115), Summer vacation (200–251) and Autumn half-term (300–307). The rate of new infections is given the form

$$(6) \qquad \mu_{SE} = \beta(t)[I(t) + \omega]^{\alpha}/P(t).$$

Here, $\omega$ describes infection from measles cases outside the population under study; $\alpha$ describes inhomogeneous mixing [Finkenstädt and Grenfell (2000)]; $P(t)$ is the total population size, which is treated as a known covariate via interpolation from census data. Environmental stochasticity on transmission is included via a gamma noise process $\xi_{SE}(t)$ with infinitesimal variance parameter $\sigma_{SE}^2$; transmission is presumed to be the most variable process in the system, and other transitions are taken to be noise-free. The two other disease parameters, $\mu_{EI}$ and $\mu_{IR}$, are treated as unknown constants. We suppose a constant mortality rate, $\mu_{SD} = \mu_{ED} = \mu_{ID} = \mu_{RD} = m$, and here we fix $m = 1/50$ year$^{-1}$. The recruitment of school-age children is specified by the process $N_{BS}(t) = \lfloor \int_0^t b(s-\tau)\,ds \rfloor$, where $b(t)$ is the birth rate, presented in Figure 3, and $\lfloor x \rfloor$ is the integer part of $x$. We note that the construction above does not perfectly match the constraint $S(t) + E(t) + I(t) + R(t) = P(t)$. For a childhood disease, such as measles, a good estimate of the birth rate is important, whereas the system is insensitive to the exact size of the adult population.

All transitions not mentioned above are taken to have a rate of zero. To complete the model specification, a measurement model is required. Biweekly aggregated measles cases are denoted by $C_n = N_{IR}(t_n) - N_{IR}(t_{n-1})$ with $t_n$ being the time, in years, of the $n$th observation. Reporting rates $\rho_n$ are taken to be independent Gamma$(1/\phi, \rho\phi)$ random variables. Conditional on $\rho_n$, the observations are modeled as independent Poisson counts, $Y_n|\rho_n, C_n \sim \text{Poisson}(\rho_n C_n)$. Thus, $Y_n$ given $C_n$ has a negative binomial distribution with $E[Y_n|C_n] = \rho C_n$ and $\text{Var}(Y_n|C_n) = \rho C_n + \phi\rho^2 C_n^2$. Note that the measurement model counts transitions into $R$, since individuals are removed from the infective pool (treated with bed-rest) once diagnosed. The measurement model allows for the possibility of both demographic stochasticity (i.e., Poisson variability) and environmental stochasticity (i.e., gamma variability on the rates).

A likelihood ratio test concludes that, in the context of this model, environmental stochasticity is clearly required to explain the data: the log likelihood for the full model was found to be $-2504.9$, compared to $-2662.0$ for the restricted model with $\sigma_{SE} = 0$ ($p < 10^{-6}$, chi-square test; results based on a time-step of $\delta = 1$ day in the Euler scheme of Figure 1 and a Monte Carlo sample size of $J = 20000$ when carrying out the iterated filtering algorithm in Figure 2). Future model-based scientific investigations of disease dynamics should consider environmental stochasticity when basing scientific conclusions on the results of formal statistical tests.

Environmental stochasticity, like over-dispersion in generalized linear models, is more readily detected than scientifically explained. Teasing apart the extent to which environmental stochasticity is describing model mis-specification rather than random phenomena in the system is beyond the scope of the present paper. The inference framework developed here will facilitate both asking and answering such questions. However, this distinction is not relevant to the central statistical question of whether a particular class of scientifically motivated models requires environmental stochasticity (in the broadest sense of a source of variability above and beyond demographic stochasticity) to explain the data. Models of biological systems are necessarily simplifications of complex processes [May (2004)], and as such, it is a legitimate role for environmental stochasticity to represent and quantify the contributions of unknown and/or unmodeled processes to the system under investigation.

The environmental stochasticity identified here has consequences for the qualitative understanding of measles epidemics. Bauch and Earn (2003) have pointed out that demographic stochasticity is not sufficient to explain the deviations which historically occurred from periodic epidemics (at one, two or three year cycles, depending on the population size and the birth-rate). Simulations from the fitted model with environmental stochasticity are able to reproduce such irregularities (results not shown), giving a simple explanation of this phenomenon. This does not rule out the possibility that some other explanation, such as explicitly introducing a new covariate into the model, could give an even better explanation.

In agreement with Cauchemez and Ferguson (2008), we have found that some combinations of parameters in our model are only weakly identifiable (i.e., they are formally identifiable, but have broad confidence intervals). Although this does not invalidate the above likelihood ratio test, it does cause difficulties interpreting parameter estimates. In the face of this problem, Cauchemez and Ferguson (2008) made additional modeling assumptions to improve identifiability of unknown parameters. Here, our goal is to demonstrate our modeling and inference framework, rather than to present a comprehensive investigation of measles dynamics.

The analysis of Cauchemez and Ferguson (2008), together with other contributions by the same authors [Cauchemez et al. (2008, 2006)], represents the state of the art for Markov chain Monte Carlo analysis of population dynamics. Whereas Cauchemez and Ferguson (2008) required model-specific approximations and analytic calculations to carry out likelihood-based inference via their Markov chain Monte Carlo approach, our analysis is a routine application of the general framework in Sections 2 and 3. Our methodology also goes beyond that of Cauchemez and Ferguson (2008) by allowing the consideration of environmental stochasticity, and the inclusion of the disease latent period (represented by the compartment $E$) which has been found relevant to the disease dynamics [Finkenstädt and Grenfell (2000), Bjornstad, Finkenstadt and Grenfell (2002), Morton and Finkenstadt (2005)]. Furthermore, our approach generalizes readily to more complex biological systems, as demonstrated by the following example.

4.2. *A mechanistic model for competing strains of cholera.* All infectious pathogens have a variety of strains, and a good understanding of the strain structure can be key to understanding the epidemiology of the disease, understanding evolution of resistance to medication, and developing effective vaccines and vaccination strategies [Grenfell et al. (2004)]. Previous analyses relating mathematical consequences of strain structure to disease data include studies of malaria [Gupta et al. (1994)], dengue [Ferguson, Anderson and Gupta (1999)], influenza [Ferguson, Galvani and Bush (2003), Koelle et al. (2006a)] and cholera [Koelle, Pascual and Yunus (2006b)]. For measles, the strain structure is considered to have negligible importance for the transmission dynamics [Conlan and Grenfell (2007)], another reason why measles epidemics form a relatively simple biological system. In this section, we demonstrate that our mechanistic modeling framework permits likelihood based inference for mechanistically motivated stochastic models of strain-structured disease systems, and that the results can lead to fresh scientific insights.

There are many possible immunological consequences of the presence of multiple strains, but it is often the case that exposure of a host to one strain of a pathogen results in some degree of protection (immunity) from re-infection by that strain and, frequently, somewhat weaker protection (cross-immunity) from infection by other strains. Immunologically distinct strains are called serotypes. In the case of cholera, there are currently two common serotypes, Inaba and Ogawa. Koelle, Pascual and Yunus (2006b), following the multistrain modeling approach of Kamo and Sasaki (2002), constructed a mechanistic, deterministic model of cholera transmission and immunity to investigate the pattern of changes in serotype dominance observed in cholera case report data collected in an intensive surveillance program conducted by the International Center for Diarrheal Disease Research, Bangladesh. They argued on the basis of a comparison of the data with features of typical trajectories of the dynamical model. Specifically, Koelle, Pascual and Yunus (2006b) found that the model would exhibit behavior which approximately matched the period of cycles in strain dominance only when the cross-immunity was high, that is, when the probability of cross-protection was approximately 0.95. In addition, they found that their model's behavior depended very sensitively on the cross-immunity parameter. Here we employ formal likelihood-based inference on the same data to assess the strength of the evidence in favor of these conclusions.

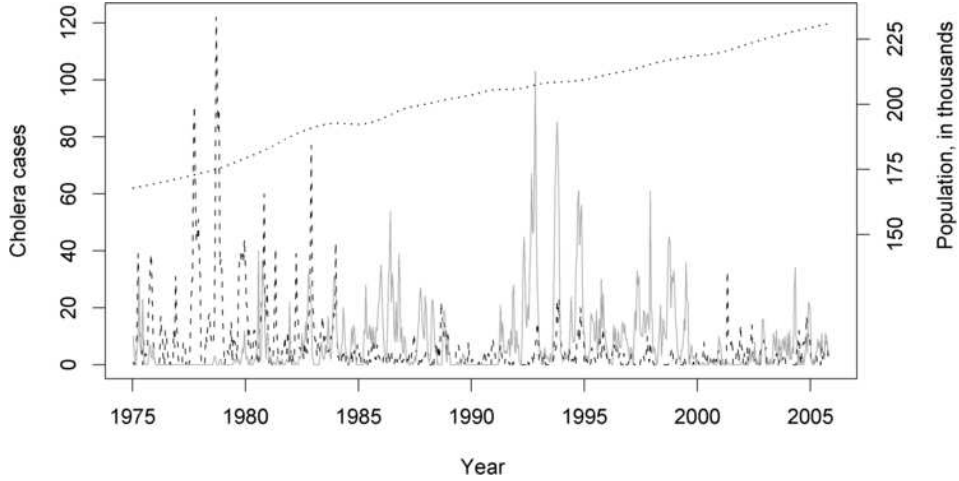

Fig. 5. *Biweekly cholera cases obtained from hospital records of the International Center for Diarrheal Disease Research, Bangladesh for the district of Matlab, Bangladesh, 55 km SE of Dhaka. Cases are categorized into serotypes, Inaba (dashed) and Ogawa (solid gray). Each serotype may be further classified into one of two biotypes, El Tor and Classical, which are combined here, following Koelle et al. (2006b). The total population size of the district, in thousands, is shown as a dotted line.*

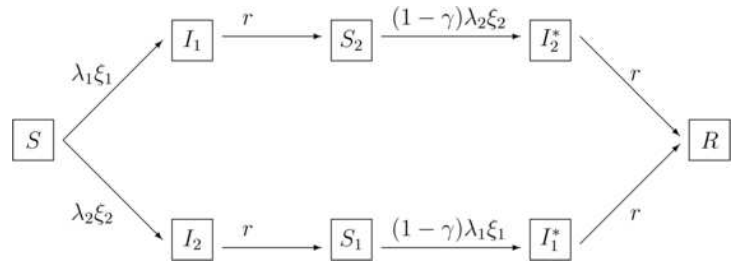

Fig. 6. *Flow diagram for cholera, including interactions between the two major serotypes. Each individual host falls in one compartment: $S$, susceptible to both Inaba and Ogawa serotypes; $I_1$, infected with Inaba; $I_2$, infected with Ogawa; $S_1$, susceptible to Inaba (but immune to Ogawa); $S_2$, susceptible to Ogawa (but immune to Inaba); $I_1^*$, infected with Inaba (but immune to Ogawa); $I_2^*$, infected with Ogawa (but immune to Inaba); $R$, immune to both serotypes. Births enter $S$, and all individuals have a mortality rate $m$.*

We analyzed a time series consisting of 30 years of biweekly cholera incidence records (Figure 5). For each cholera case, the serotype of the infecting strain—Inaba or Ogawa—was determined. We formulated a stochastic version of the model analyzed by Koelle, Pascual and Yunus (2006b). The model is shown diagrammatically in Figure 6, in which arrows represent possible transitions, each labeled with the corresponding rate of flow. Table 1 specifies the model formally, in the framework of Section 2, as a Markov chain with stochastic rates. The parameters in the model have standard epidemiological interpretations [Anderson and May (1991), Finkenstädt and Grenfell (2000), Koelle and Pascual (2004)]: $\lambda_1$ is the force of infection for the Inaba serotype, that is, the mean rate at which susceptible individuals become infected; $\xi_1$ is the stochastic noise on this rate; $\lambda_2$ and $\xi_2$ are the corresponding force of infection and noise for the Ogawa serotype; $\beta(t)$ is the rate of transmission between individuals, parameterized with a linear trend and a smooth seasonal component; $\omega$ gives the rate of infection from an environmental reservoir, independent of the current number of contagious individuals; the exponent $\alpha$ allows for inhomogeneous mixing of the population; $r$ is the recovery rate from infection; $\gamma$ measures the strength of cross-immunity between serotypes. In this model, as in Koelle, Pascual and Yunus (2006b), infection with a given serotype results in life-long immunity to reinfection by that serotype. The argument for giving both strains common variability is that they are believed to be biologically similar except in regard to immune response. The strains have independent noise components because the noise represents chance events, such as a contaminated feast or a single community water source which is transiently in a favorable condition for contamination, and such events spread whichever strain is in the required place at the required time.

To complete the model specification, we adopt an extension of the negative binomial measurement model used for measles in Section 4.1. Biweekly aggregated cases for Inaba and Ogawa strains are denoted by $C_{i,n} = N_{SI_i}(t_n) - N_{SI_i}(t_{n-1}) + N_{S_iI_i^*}(t_n) - N_{S_iI_i^*}(t_{n-1})$ for $i = 1, 2$ respectively and $t_n = 1975 + n/24$. Reporting rates $\rho_{1,n}$ and $\rho_{2,n}$ are taken to be independent Gamma$(1/\phi, \rho\phi)$ random variables. Conditional on $\rho_{1,n}$ and $\rho_{2,n}$, the observations are modeled as independent Poisson counts,

$$Y_{i,n}|\rho_{i,n}, C_{i,n} \sim \text{Poisson}(\rho_{i,n}C_{i,n}), \qquad i = 1, 2.$$

Thus, $Y_{i,n}$ given $C_{i,n}$ has a negative binomial distribution with $E[Y_{i,n}|C_{i,n}] = \rho C_{i,n}$ and $\text{Var}(Y_{i,n}|C_{i,n}) = \rho C_{i,n} + \phi\rho^2 C_{i,n}^2$.

Some results from fitting the model in Figure 6 via the method in Figure 2 are shown in Table 2. The two sets of parameter values $\hat\theta_A$ and $\hat\theta_B$ in Table 2 are maximum likelihood estimates, with $\hat\theta_A$ having the additional constraints $\rho = 0.067$ and $r = 38.4$. These two constraints were imposed by

<div style="text-align:center">

TABLE 1

*Interpretation of Figure 6 via the multinomial process with random rates in (2), with*
$X(t) = (S(t),\ I_1(t),\ I_2(t),\ S_1(t),\ S_2(t),\ I_1^*(t),\ I_2^*(t),\ R(t),\ B(t),\ D(t))$. *Compartments B and D are introduced for demographic considerations: births are formally treated as transitions from B to S and deaths as transitions into D. All transitions not listed above have zero rate.* $\xi_2(t)$ *and* $\xi_1(t)$ *are independent gamma noise processes, each with infinitesimal variance parameter* $\sigma^2$. *Transition rates are noise-free unless specified otherwise. Seasonality is modeled via a periodic cubic B-spline basis* $\{s_i(t), i = 1,\ldots,6\}$, *where* $s_i(t)$ *attains its maximum at* $t = (i-1)/6$. *The population size* $P(t)$ *is shown in Figure 5. The birth process is treated as a covariate, that is, the analysis is carried out conditional on the process* $N_{BS} = \lfloor P(t) - P(0) + \int_0^t mP(s)\,ds \rfloor$, *where* $\lfloor x \rfloor$ *is the integer part of x. There is a small stochastic discrepancy between* $S(t) + I_1(t) + I_2(t) + S_1(t) + S_2(t) + I_1^*(t) + I_2^*(t) + R(t)$ *and* $P(t)$. *In principle, one could condition on the demographic data by including a population measurement model—we saw no compelling reason to add this extra complexity for the current purposes. Numerical solutions of sample paths were calculated using the algorithm in Figure 1, with* $\delta = 2/365$

</div>

$$\lambda_1 = \beta(t)(I_1(t) + I_1^*(t))^\alpha / P(t) + \omega$$
$$\lambda_2 = \beta(t)(I_2(t) + I_2^*(t))^\alpha / P(t) + \omega$$
$$\log \beta(t) = b_0(t - 1990) + \sum_{i=1}^{6} b_i s_i(t)$$
$$\mu_{SI_1} = \lambda_1$$
$$\mu_{SI_2} = \lambda_2$$
$$\mu_{S_1 I_1^*} = (1 - \gamma)\lambda_1$$
$$\mu_{S_2 I_2^*} = (1 - \gamma)\lambda_2$$
$$\mu_{I_1 S_1} = \mu_{I_2 S_1} = r$$
$$\mu_{I_2^* R} = \mu_{I_1^* R} = r$$
$$\mu_{X_j D} = m \quad \text{for } X_j \in \{S, I_1, I_2, S_1, S_2, I_1^*, I_2^*, R\}$$
$$\xi_{SI_2} = \xi_{S_2 I_2^*} = \xi_2(t)$$
$$\xi_{SI_1} = \xi_{S_1 I_1^*} = \xi_1(t)$$

Koelle, Pascual and Yunus (2006b), based on previous literature. The fitted model with the additional constraints is qualitatively different from the unconstrained model, and we refer to the neighborhoods of these two parameter sets as regimes $A$ and $B$. Figure 7 shows a profile likelihood for cross-immunity in regime $A$; this likelihood-based analysis leads to substantially lower cross-immunity than the estimate of Koelle, Pascual and Yunus (2006b). In regime $B$, the cross-immunity is estimated as being complete ($\gamma = 1$), however, the corresponding standard error is large: cross-immunity is poorly identified in regime $B$ since the much higher reporting rate ($\rho = 0.65$) means that there are many fewer cases and so few individuals are ever exposed to both serotypes.

These two regimes demonstrate two distinct uses of a statistical model—first, to investigate the consequences of a set of assumptions and, second, to challenge those assumptions. If we take for granted the published estimates of certain parameter values, the resulting parameter estimates $\hat\theta_A$ are broadly consistent with previous models for cholera dynamics in terms of

*Parameter estimates from both regimes. In both regimes, the mortality rate m is fixed at $1/38.8$ years$^{-1}$. The units of $r$, $b_0$ and $\omega$ are year$^{-1}$; $\sigma$ has units year$^{1/2}$; and $\rho$, $\gamma$, $\phi$, $\alpha$ and $b_1, \ldots, b_6$ are dimensionless. $\ell$ is the average of two log-likelihood evaluations using a particle filter with 120,000 particles. Optimization was carried out using the iterated filtering in Figure 2, with $M = 30$, $a = 0.95$ and $J = 15{,}000$. Optimization parameters were selected via diagnostic convergence plots [Ionides, Bretó and King (2006)]. Standard errors (SEs) were derived via a Hessian approximation; this is relatively rapid to compute and gives a reasonable idea of the scale of uncertainty, but profile likelihood based confidence intervals are more appropriate for formal inference*

|  | $\hat{\theta}_A$ | $SE_A$ | $\hat{\theta}_B$ | $SE_B$ |
|---|---|---|---|---|
| $r$ | 38.42 | – | 36.91 | 3.88 |
| $\rho$ | 0.067 | – | 0.653 | 0.069 |
| $\gamma$ | 0.400 | 0.087 | 1.00 | 0.41 |
| $\sigma$ | 0.1057 | 0.0076 | 0.0592 | 0.0075 |
| $\phi$ | 0.014 | 0.30 | 0.0004 | 0.024 |
| $\omega \times 10^3$ | 0.099 | 0.022 | 0.0762 | 0.0072 |
| $\alpha$ | 0.860 | 0.015 | 0.864 | 0.017 |
| $b_0$ | $-0.0275$ | 0.0017 | $-0.0209$ | 0.0015 |
| $b_1$ | 4.608 | 0.098 | 3.507 | 0.083 |
| $b_2$ | 5.342 | 0.074 | 3.733 | 0.091 |
| $b_3$ | 5.723 | 0.075 | 4.448 | 0.055 |
| $b_4$ | 5.022 | 0.076 | 3.534 | 0.065 |
| $b_5$ | 5.508 | 0.064 | 4.339 | 0.053 |
| $b_6$ | 5.804 | 0.059 | 4.274 | 0.039 |
| $\ell$ | $-3560.23$ | | $-3539.11$ | |

the fraction of the population which has acquired immunity to cholera and the fraction of cases which are asymptomatic (a term applied to the majority of unreported cases which are presumed to have negligible symptoms but are nevertheless infectious). However, the likelihood values in Table 2 call into question the assumptions behind regime $A$, since the data are better explained by regime $B$ ($p < 10^{-6}$, likelihood ratio test), for which the epidemiologically relevant cases are only the severe cases that are likely to result in hospitalization. Unlike in regime $A$, asymptomatic cholera cases play almost no role in regime $B$, since the reporting rate is an order of magnitude higher. The contrast between these regimes highlights a conceptual limitation of compartment models: in point of fact, disease severity and level of infectiousness are continuous, not discrete or binary as they must be in basic compartment models. For example, differences in the level of morbidity required to be classified as "infected" result in re-interpretation of the parameters of the model, with consequent changes to fundamental model characteristics such as the basic reproductive ratio of an infectious disease [Anderson and May (1991)]. Despite this limitation, it remains the case that compartment models are a fundamental platform for current understanding of disease dynamics and so identification and comparison of different interpretations is an important exercise.

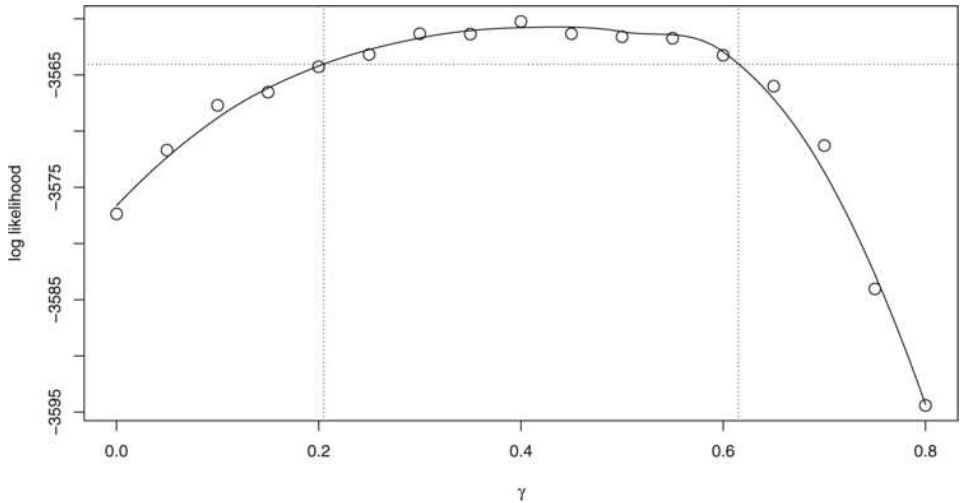

FIG. 7.   *Cross-immunity profile likelihood for regime A, yielding a 99% confidence interval for γ of (0.20, 0.61) based on a χ² approximation [e.g., Barndorff-Nielsen and Cox (1994)]. Each point corresponds to an optimization carried out as described in the caption of Table 2. Local quadratic regression implemented in R via loess with a span of 0.6 [R Development Core Team (2006)] was used to estimate the profile likelihood, following Ionides (2005).*

The existence of regime B shows that there is room for improvement in the model by departing from the assumptions in regime A. Extending the model to include differing levels of severity might permit a combination of the data-matching properties of B with the scientific interpretation of A. Other aspects of cholera epidemiology [Sack et al. (2004), Kaper, Morris and Levine (1995)], not included in the models considered here, might affect the conclusions. For example, although we have followed Koelle, Pascual and Yunus (2006b) by assuming lifelong immunity following exposure to cholera, in reality this protection is believed to wane over time [King et al. (2008b)]. It is plain, however, that any such model modification can be subsumed in the modeling framework presented here and that effective inference for such models is possible using the same techniques.

**5. Discussion.** This paper has focused on compartment models, a flexible class of models which provides a broad perspective on the general topic of mechanistic models. The developments of this paper are also relevant to other systems. For example, compartment models are closely related to chemodynamic models, in which a Markov process is used to represent the quantities of several chemical species undergoing transformations by chemical reactions. The discrete nature of molecular counts can play a role, particularly for large biological molecules [Reinker, Altman and Timmer (2006), Boys, Wilkinson and Kirkwood (2008)]. Our approach to stochastic transition rates (Section 2) could readily be extended to chemodynamic models, and would allow for the possibility of over-dispersion in experimental systems.

Given rates $\mu_{ij}$, one interpretation of a compartment model is to write the flows as coupled ordinary differential equations (ODEs),

$$\frac{d}{dt} N_{ij} = \mu_{ij} X_i(t). \tag{7}$$

Data analysis via ODE models has challenges in its own right [Ramsay et al. (2007)]. One can include stochasticity in (7) by adding a slowly varying function to the derivative [Swishchuk and Wu (2003)]. Alternatively, one can add Gaussian white noise to give a set of coupled stochastic differential equations (SDEs) [e.g., Øksendal (1998)]. For example, if $\{W_{ijk}(t)\}$ is a collection of independent standard Brownian motion processes, and $\sigma_{ijk} = \sigma_{ijk}(t, X(t))$, an SDE interpretation of a compartment model is given by

$$dN_{ij} = \mu_{ij} X_i(t)\,dt + \sum_k \sigma_{ijk}\,dW_{ijk}.$$

SDEs have some favorable properties for mechanistic modeling, such as the ease with which stochastic models can be written down and interpreted in terms of infinitesimal mean and variance [Ionides et al. (2004), Ionides, Bretó and King (2006)]. However, there are several reasons to prefer integer-valued stochastic processes over SDEs for modeling population processes. Populations consist of discrete individuals, and, when a population becomes small, that discreteness can become important. For infectious diseases, there may be temporary extinctions, or "fade-outs," of the disease in a population or sub-population. Even if the SDE is an acceptable approximation to the disease dynamics, there are technical reasons to prefer a discrete model. Standard methods allow exact simulation for continuous time Markov chains [Brémaud (1999), Gillespie (1977)], whereas for a nonlinear SDE this is at best difficult [Beskos et al. (2006)]. In addition, if an approximate Euler solution for a compartment model is required, nonnegativity constraints can more readily be accommodated for Markov chain models, particularly when the model is specified by a limit of multinomial approximations, as in (2). The most basic discrete population compartment model is the Poisson system [Brémaud (1999)], defined here by

$$P[\Delta N_{ij} = n_{ij} | X(t) = (x_1, \ldots, x_c)]$$
$$= \prod_i \prod_{j \neq i} (\mu_{ij} x_i \delta)^{n_{ij}} (1 - \mu_{ij} x_i \delta) + o(\delta). \tag{8}$$

The Poisson system is a Markov chain whose transitions consist of single individuals moving between compartments, that is, the infinitesimal probability is negligible of either simultaneous transitions between different pairs of compartments or multiple transitions between a given pair of compartments. As a consequence of this, the Poisson system is "equidispersed," meaning that the infinitesimal mean of the increments equals the infinitesimal variance (Appendix B). Overdispersion is routinely observed in data [McCullagh and Nelder (1989)], and this leads us to consider models such as (2) for which the infinitesimal variance can exceed the infinitesimal mean. For infinitesimally over-dispersed systems, instantaneous transitions of more than one individual are possible. This may be scientifically plausible: a cholera-infected meal or water-jug may lead to several essentially simultaneous cases; many people could be simultaneously exposed to an influenza patient on a crowded bus. Quite aside from this, if one wishes for whatever reason to write down an over-dispersed Markov model, the inclusion of such possibilities is unavoidable. Simultaneity in the limiting continuous time model can alternatively be justified by arguing that the model only claims to capture macroscopic behavior over sufficiently long time intervals.

Note that the multinomial distribution used in (2) could be replaced by alternatives, such as Poisson or negative binomial. These alternatives are more natural for unbounded processes, such as birth processes. For equidispersed processes, that is, without adding white noise to the rates, the limit in (2) is the same if the multinomial is replaced by Poisson or negative binomial. For overdispersed processes, these limits differ. In particular, the Poisson gamma and negative binomial gamma limits have unbounded jump distributions and so are less readily applicable to finite populations.

The approach in (2) of adding white noise to the transition rates differs from previous approaches of making the rates a slowly varying random function of time, that is, adding low frequency "red noise" to the rates. There are several motivations for introducing models based on white noise. Most simply, adding white noise can lead to more parsimonious parameterization, since the intensity but not the spectral shape of the noise needs to be considered. The Markov property of white noise is inherited by the dynamical system, allowing the application of the extensive theory of Markov chains. White noise can also be used as a building block for constructing colored noise, for example, by employing an autoregressive model for the parameters. At least for the specific examples of measles and cholera studied in Section 4, high-frequency variability in the rate of infection helps to explain the data (this was explicitly tested for measles; for our cholera models, Table 2 shows that the estimates of the environmental stochasticity parameter $\sigma$ are many standard errors from zero). Although variability in rates will not always be best modeled using white noise, there are many circumstances in which it is useful to be able to do so.

This article has taken a likelihood-based, non-Bayesian approach to statistical inference. Many of the references cited follow the Bayesian paradigm. The examples of Section 4 and other recent work [Ionides, Bretó and King (2006), King et al. (2008b)] show that iterated filtering methods enable routine likelihood-based inference in some situations that have been challenging for Bayesian methodology. Bayesian and non-Bayesian analyses will continue to provide complementary approaches to inference for time series analysis via mechanistic models, as in other areas of statistics.

Time series analysis is, by tradition, data oriented, and so the quantity and quality of available data may limit the questions that the data can reasonably answer. This forces a limit on the number of parameters that can be estimated for a model. Thus, a time series model termed mechanistic might be a simplification of a more complex model which more fully describes reductionist scientific understanding of the dynamical system. As one example, one could certainly argue for including age structure or other population inhomogeneities into Figure 6. Indeed, determining which additional model components lead to important improvement in the statistical description of the observed process is a key data analysis issue.

## APPENDIX A: THEOREMS CONCERNING COMPARTMENT MODELS WITH STOCHASTIC RATES

THEOREM A.1. *Assume* (P1)–(P5) *and suppose that* $\mu_{ij}(t, x)$ *is uniformly continuous as a function of* $t$. *Suppose initial values* $X(0) = (X_1(0), \ldots, X_c(0))$ *are given and denote the total number of individuals in the population by* $S = \sum_i X_i(0)$. *Label the individuals* $1, \ldots, S$ *and the compartments* $1, \ldots, c$. *Let* $C(\zeta, 0)$ *be the compartment containing individual* $\zeta$ *at time* $t = 0$. *Set* $\tau_{\zeta,0} = 0$, *and generate independent* Exponential(1) *random variables* $M_{\zeta,0,j}$ *for each* $\zeta$ *and* $j \neq C(\zeta, 0)$. *We will define* $\tau_{\zeta,m,j}$, $\tau_{\zeta,m}$, $C(\zeta, m)$, *and* $M_{\zeta,m,j}$ *recursively for* $m \geq 1$. *Set*

$$(9) \quad \tau_{\zeta,m,j} = \inf\left\{t : \int_{\tau_{\zeta,m-1}}^{t} \mu_{C(\zeta,m-1),j}(s, X(s)) \, d\Gamma_{C(\zeta,m-1),j}(s) > M_{\zeta,m-1,j}\right\}.$$

*At time* $\tau_{\zeta,m} = \min_j \tau_{\zeta,m,j}$, *set* $C(\zeta, m) = \arg\min_j \tau_{\zeta,m,j}$ *and for each* $j \neq C(\zeta, m)$ *generate an independent* Exponential(1) *random variable* $M_{\zeta,m,j}$. *Set*

$$dN_{ij}(t) = \sum_{\zeta,m} \mathbb{I}\{C(\zeta, m-1) = i, C(\zeta, m) = j, \tau_{\zeta,m} = t\},$$

*where* $\mathbb{I}\{\cdot\}$ *is an indicator function, and set* $X_i(t) = X_i(0) + \int_0^t \sum_{j \neq i}(dN_{ji} - dN_{ij})$. $X(t)$ *is a Markov chain whose infinitesimal transition probabilities satisfy* (2).

NOTE A.1. The random variables $\{M_{\zeta,m,j}\}$ in Theorem A.1 are termed *transition clocks*, with the intuition that $X(t)$ jumps when one (or more) of the integrated transition rates in (9) exceeds the value of its clock. In a more basic construction of a Markov chain, one re-starts the clocks for each individual whenever $X(t)$ makes a transition. The memoryless property of the exponential distribution makes this equivalent to the construction of Theorem A.1, where clocks are restarted only for individuals who make a transition [Sellke (1983)]. Sellke's construction is convenient for the proof of Theorem A.1.

NOTE A.2. The trajectories of the individuals are coupled through the dependence of $\mu(t, X(t))$ on $X(t)$, and through the noise processes which are shared for all individuals in a given compartment and may be dependent between compartments. Thus, to evaluate (9), it is necessary to keep track of all individuals simultaneously. To check that the integral in (9) is well defined, we note that $X(s)$ depends only on $\{(\tau_{\zeta,m}, C(\zeta, m)) : \tau_{\zeta,m} \leq s, \zeta = 1, \ldots, S, m = 1, 2, \ldots\}$. $X(s)$ is thus a function of events occurring by time $s$, so it is legitimate to use $X(s)$ when constructing events that occur at $t > s$.

THEOREM A.2. *Supposing* (P1)–(P7), *the infinitesimal transition probabilities given by* (2) *are*

$$P[\Delta N_{ij} = n_{ij}, for\ all\ i \neq j \mid X(t) = (x_1, \ldots, x_c)]$$
$$= \prod_i \prod_{j \neq i} \pi(n_{ij}, x_i, \mu_{ij}, \sigma_{ij}) + o(\delta),$$

*where*

(10)
$$\pi(n, x, \mu, \sigma)$$
$$= 1_{\{n=0\}} + \delta \binom{x}{n} \sum_{k=0}^{n} \binom{n}{k} (-1)^{n-k+1} \sigma^{-2} \ln(1 + \sigma^2 \mu(x - k)).$$

The full independence of $\{\Gamma_{ij}\}$ assumed in Theorem A.2 gives a form for the limiting probabilities where multiple individuals can move simultaneously between some pair of compartments $i$ and $j$, but no simultaneous transitions occur between different compartments. In more generality, the limiting probabilities do not have this simple structure. In the setup for Theorem A.1, where $\Gamma_{ij}$ is independent of $\Gamma_{ik}$ for $j \neq k$, no simultaneous transitions occur out of some compartment $i$ into different compartments $j \neq k$, but simultaneous transitions from $i$ to $j$ and from $i'$ to $j'$ cannot be ruled out for $i \neq i'$. The assumption in Theorem A.1 that $\Gamma_{ij}$ is independent of $\Gamma_{ik}$ for $j \neq k$ is not necessary for the construction of a process via (2), but simplifies the subsequent analysis. Without this assumption, a construction similar to Theorem A.1 would have to specify a rule for what happens when an individual who has two simultaneous event times, that is, when $\min_j \tau_{\zeta,m,j}$ is not uniquely attained. Although independence assumptions are useful for analytical results, a major purpose of the formulation in (2) is to allow the practical use of models that surpass currently available mathematical analysis. In particular, it may be natural for different transition processes to share the same noise process, if they correspond to transitions between similar pairs of states.

## APPENDIX B: EQUIDISPERSION OF POISSON SYSTEMS

For the Poisson system in (8),

$$P[\Delta N_{ij} = 0 | X(t) = x] = 1 - \mu_{ij} x_i \delta + o(\delta), \tag{11}$$

$$P[\Delta N_{ij} = 1 | X(t) = x] = \mu_{ij} x_i \delta + o(\delta), \tag{12}$$

where $x = (x_1, \ldots, x_c)$ and $\mu_{ij} = \mu_{ij}(t, x)$. Since the state space of $X(t)$ is finite, it is not a major restriction to suppose that there is some uniform bound $\mu_{ij}(t, x) x_i \leq \nu$, and that the terms $o(\delta)$ in (11), (12) are uniform in $x$ and $t$. Then, $P[\Delta N_{ij} > k | X(t)] \leq \bar{F}(k, \delta \nu)$, where $\bar{F}(k, \lambda) = \sum_{j=k+1}^{\infty} \lambda^j e^{-\lambda}/j!$. It follows that

$$E[\Delta N_{ij} | X(t) = x] = \sum_{k=0}^{\infty} P[\Delta N_{ij} > k | X(t) = x] = \mu_{ij} x_i \delta + o(\delta),$$

$$E[(\Delta N_{ij})^2 | X(t) = x] = \sum_{k=0}^{\infty} (2k+1) P[\Delta N_{ij} > k | X(t) = x] = \mu_{ij} x_i \delta + o(\delta),$$

and so $\mathrm{Var}(\Delta N_{ij} | X(t) = x) = \mu_{ij} x_i \delta + o(\delta)$. If the rate functions $\mu_{ij}(t, x)$ are themselves stochastic, with $X(t)$ being a conditional Markov chain given $\{\mu_{ij}, 1 \leq i \leq c, 1 \leq j \leq c\}$, a similar calculation applies so long as a uniform bound $\nu$ still exists. In this case,

$$E[\Delta N_{ij} | X(t) = x] = \delta E[\mu_{ij}(t, x) x_i] + o(\delta), \tag{13}$$

$$\mathrm{Var}(\Delta N_{ij} | X(t) = x) = \delta E[\mu_{ij}(t, x) x_i] + o(\delta). \tag{14}$$

The necessity of the uniform bound $\nu$ is demonstrated by the inconsistency between (13), (14) and the result in Theorem A.2 for the addition of white noise to the rates.

**Acknowledgments.** This work was stimulated by the working groups on "Inference for Mechanistic Model" and "Seasonality of Infectious Diseases" at the National Center for Ecological Analysis and Synthesis, a Center funded by NSF (Grant DEB-0553768), the University of California, Santa

Barbara, and the State of California. We thank Mercedes Pascual for discussions on multi-strain modeling. We thank the editor, associate editor and referee for their helpful suggestions. The cholera data were kindly provided by the International Centre for Diarrheal Disease Research, Bangladesh.

## SUPPLEMENTARY MATERIAL

**Theorems concerning compartment models with stochastic rates** (DOI: 10.1214/08-AOAS201SUPP; .pdf). We present proofs of Theorems A.1 and A.2, which were stated in Appendix A.

C. BRETÓ
DEPARTMENT OF STATISTICS
UNIVERSIDAD CARLOS III DE MADRID
GETAFE 28903
MADRID
SPAIN
E-MAIL: cbreto@est-econ.uc3m.es

D. HE
A. A. KING
DEPARTMENT OF ECOLOGY
  AND EVOLUTIONARY BIOLOGY
UNIVERSITY OF MICHIGAN
ANN ARBOR, MICHIGAN 48109-1048
USA
E-MAIL: daihai@umich.edu
        aaron.king@umich.edu

E. L. IONIDES
DEPARTMENT OF STATISTICS
UNIVERSITY OF MICHIGAN
ANN ARBOR, MICHIGAN 48109-1107
USA
E-MAIL: ionides@umich.edu